\documentclass{article}
\usepackage[utf8]{inputenc}
\usepackage{tabularx}
\usepackage{amsmath}
\usepackage{amssymb}
\usepackage{amsthm}
\usepackage{geometry}
\usepackage{enumitem}
\usepackage[T1]{fontenc}
\usepackage{graphicx}
\usepackage[dvipsnames]{xcolor}
\title{On the eleventh degree transformation of elliptic functions }
\author{Felix Klein, Translation by Yonathan Stone\thanks{The translator would like to thank Jesse Wolfson for suggesting this and for providing helpful feedback on earlier drafts. In addition, this translation was supported in part by NSF Grant DMS-1944862.}}
\date{Mathematische Annalen, Vol 15 (1879), Translated in 2020}

\begin{document}
\maketitle
In the pursuit of my study of the transformation of elliptic functions I will treat the case of the \textit{eleventh} degree transformation in the sequel.  In doing so it bears mentioning that I wish to explicitly produce the \textit{eleventh} degree equation, which arises in this case, in its simplest form. In the XIV\textsuperscript{th} volume of these annals, pp. 423-424, I have already shown that this equation can be presented in the following form:
\[J = F(z),\]
where $F(z)$ is an \textit{entire} function in $z$ of eleventh degree with only numerical coefficients, which contains a cubic factor three times, while $(F(z) - 1)$ possesses twofold biquadratic factor.  But at the same time I must remark that these specifications does not suffice for completely determining the function $F$.  Thus I will combine the aforementioned results with a result from the XV\textsuperscript{th} volume pg. 277 concerning transformations of any degree.  \textit{The same with $\frac{n-1}{2}$ variables $y$ gives rise to a system of collineations that is isomorphic to the group of the modular equation.}  Corresponding to this, we have a ``Problem of $y$'' of 660\textsuperscript{th} degree. A suitable specialization of this firstly gives us the Galois resolvent of 660\textsuperscript{th} degree of the modular equation in a very straightforward manner, then subsequently the desired equation of eleventh degree. This takes on two forms, each of which has its respective advantages, where one is $J = F(z)$, as above.  $\S10$ contains my summary of these results.  The following sections communicate the transition to the simple multiplier equation of \textit{twelfth} degree, which I outlined in Vol. XV, pg. 88, and provide the possibility of solving the new equations of eleventh and 660\textsuperscript{th} degrees in a transcendental manner.\\
I've already publicized the results listed herein, albeit without proof, in two letters to Mr. Brioschi on the 9 April and 11 June 1879 in the \textit{Academia dei Lincei} and \textit{Istituto Lombardo} respectively.  On the other hand the entire development as given below, was given in a lecture on algebraic equations during the previous summer semester.
\newpage
\textbf{\S1 On certain eleven sheeted Riemann surfaces\footnote{This section is linked unambiguously to the research referenced previously: \textit{Ueber die Erniedrigung der Modulargleichungen}, Ann. XIV, pp. 417-427. [More general computations of the number of Riemann surfaces possessing prescribed branch points in the $z$-plane can be found primarily in the two papers due to Hurwitz in the 39\textsuperscript{th} and 55\textsuperscript{th} volume of the Annalen (1891 and 1901/02). K.]}}\\
The root $z$ of the equation
\begin{equation}
    J = F(z) \tag{1}
\end{equation}
which was just discussed is, as I ultimately showed, branched in relation to $J$ in such a way that all eleven sheets are associated cyclically at $J = \infty$, three times three sheets at $J = 0$, and four times two sheets at $J = 1$.  I claim at the same time, \textit{that there exist no less than ten essentially different Riemann surfaces possessing the aforementioned property} (of which only two are considered in the theory of transformations).  Now, my next task is to prove this claim in the purely geometric fashion which was previously alluded to.\\
As I did previously for the sake of clarity, I would like to shade half of each of the eleven sheets of the Riemann surface, namely to the extent that they cover the positive half plane of $J$.  Next one cuts each of the sheets along the real axis of $J$, which ranges from $J = 0$ to $J = 1$ in the finite part.  Following from the prescribed multiplicity of the branching points, our surface thus becomes a simply connected, simply bounded surface, whose sheets are still connected in a cycle at $J = \infty$.  Clearly one can use a smooth deformation to stretch this out to the interior of a circle that is segmented into 22 alternating shaded and unshaded sectors from any one of its points.  If one wants to arrange the points $J = 0$ as corners and to draw the circle arcs between them as straight lines, one has the 11-gon found in figure (8) of the panel found accompanying my previous work (Vol. XIV).
\\
Solely in consideration of the special purpose of the manner in which the newly included Fig. 1. is meant to help visualize, I will modify the figure that was just created.  Rather than considering the \textit{interior} of the 11-gon as the image of the cut up surface, I will choose the \textit{exterior} 
\\
\includegraphics[width=\textwidth]{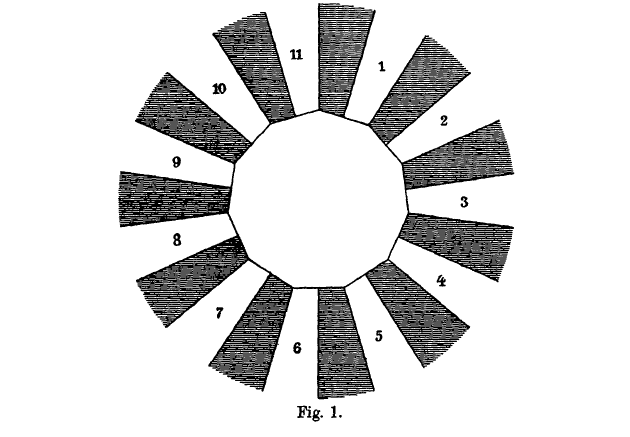} 
\\
of the same and replace the 22 half diagonals, previously converging at the center, by just as many straight lines moving towards infinity.  Besides this, for the sake of designating the different sheets, I have labeled the non-shaded regions with the numbers 1,2,...,11.\\
If one wishes to subsequently know how many different eleven-sheeted surfaces there are whose branch points have the locations and multiplicity we prescribed, then the question is clearly the following (see p. 86 of the present volume): \textit{In how many different ways is it possible to take the 22 half-edges of the inner boundary in Fig. 1 and bend them together into a doubly covered polygonal chain consisting of 11 pieces, such that of the 11 points three come together three times at $J = 0$, two come together four times at $J = 1$?} Fig. 3 on page 144 (which is surely comprehensible without further explanation) should illustrate using an example what is meant by this process of bending together.\\
The question formulated in this manner is addressed by the diagrams I,...,VI in Figure 2.  The shape of the polygonal path of length eleven should correspond only to these in every case.  The points marked with small circles each correspond to $J = 0$, and those marked by straight strokes to $J = 1$.  Diagram V relates to the case shown in Fig. 3, the diagrams I are, following from my previous explanations, the only ones which fit the equations of eleventh degree stemming from transformation theory.\\
\begin{center}
\includegraphics[width = 8.5cm]{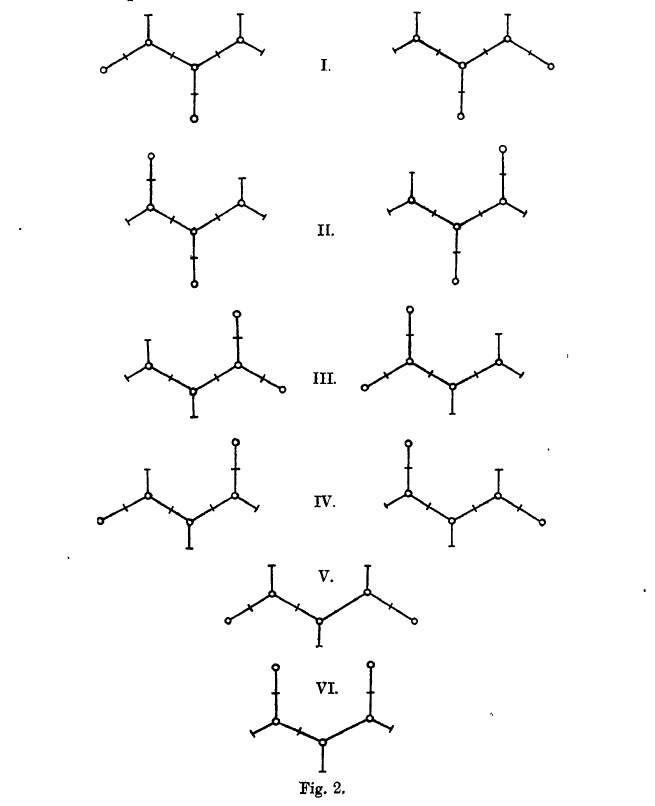}
\end{center}

According to these diagrams there are indeed \textit{ten} different ways to do the bending procedure.  The fact that there aren't also more is similarly clear; obviously it is impossible, to create new polygonal paths of length eleven with the desired properties.  Therefore, the theorem stated at the start of this section has been proven true.
\\
\textbf{\S2 Monodromy Group}
\\
As from now, we will deal with finding a simple criterion for characterizing the only two cases out the ten which bear importance for transformation theory.  I have chosen this to be \textit{Monodromy group}\footnote{See footnote \textsuperscript{22} on page 135 of the present volume}.  If one considers arbitrary closed paths in the plane of $J$, one knows that in the two cases being considered the 11 roots of the equation $J = F(z)$ are permuted in only 660 different ways, and the 660 permutations form the known group, initially discovered by Galois.  Betti\footnote{Annali di Scienze matematiche etc. di Tortolini, t. IV (1853) [= Opere matematiche, No. VI, part I., p. 81ff]} would later publish the first comprehensive investigation of this group.  This group contains only those permutations whose period is 1,2,3,5,6, and 11; and I want to use this fact here to show that the monodromy group differs from this one for the eight cases unusable to us.  \textit{That is, both surfaces of interest are completely characterized by their branch points and monodromy group.}
\\
The proof of this is analogous in all cases.  I will thus only clarify the case of Fig. 3 (diagram V).  Let $J$ 
\\
\includegraphics[width = \textwidth]{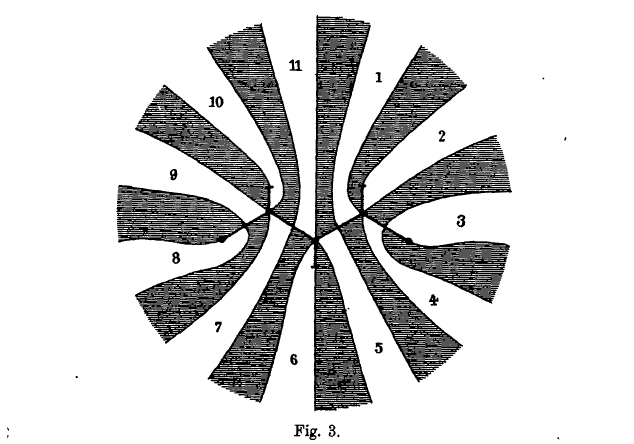}
\\
circle the point at infinity in its plane once, an operation I will denote by $S$.  As exhibited in the figure, the roots
\[1,2,3,4,5,6,7,8,9,10,11,\]
become (under suitable choice of orientation) the following:
\[2,3,4,5,6,7,8,9,10,11,1.\]
On the other hand, let $J$ circle the point $J = 1$, and call this operation $T$.  Thus,
\[1,2,3,4,5,6,7,8,9,10,11,\]
respectively become, as the figure shows,
\[5,2,4,3,1,6,11,9,8,10,7.\]
Now, we will first apply the operation $T$ once, followed by the operation $S$ three times.  Starting with
\[1,2,3,4,5,6,7,8,9,10,11,\]
we end with 
\[8,5,7,6,4,9,3,1,11,2,10;\]
\textit{Therefore, the operation $TS^3$ permutes the roots according to the following cycle decomposition:}
\[(1,8)(2,5,4,6,9,11,10)(3,7).\]
The cycles contain either 2 or 7 symbols;  The period of $TS^3$ is thus 14, and the group of monodromy thus differs from the group of 660 substitutions, which is what was to be proven\footnote{In the other cases it also suffices to examine the operation $TS^3$.}.
\\
\textbf{\S3 The Group of 660 y-substitutions.}
\\
I will utilize this result by way of considering a problem which always possesses the right group of 660 permutations; \textit{if in the process of this investigation I succeed at arriving at an equation $J = F(z)$, where $F(z)$ contains a cubic factor threefold and $F(z) - 1$ contains a biquadratic factor twice, I will have found the equation being sought after.}
\\
For this purpose I will consider the system of 660 collineations in 5 variables, as was mentioned in the introduction.  By choosing the quadratic residues modulo 11 as indices and ordering them by how they arise from multiplication by 4, I will call the variables
\[y_1,y_4,y_5,y_9,y_3.\]
The 660 collineations then emerge by repeating and combining the following two operations ($\rho = e^{\frac{2\pi i}{11}}$):
\begin{equation*}
    \begin{cases}S) \quad y_1' = \rho y_1, y_4' = \rho^4y_4, y_5' = \rho^5y_5, y_9' = \rho^9y_9, y_3' = \rho^3y_3 \\ 
    T) \begin{cases}\sqrt{-11}\cdot y_1' = (\rho^9 - \rho^2)y_1 + (\rho^4 - \rho^7)y_4 + (\rho^3 - \rho^8)y_5 + (\rho^5 - \rho^6)y_9 + (\rho^1 - \rho^{10})y_3,\\ \sqrt{-11}\cdot y_4' = (\rho^4 - \rho^7)y_1 + (\rho^3 - \rho^8)y_4 + (\rho^5 - \rho^6)y_5 + (\rho^1 - \rho^{10})y_9 + (\rho^9 - \rho^{2})y_3,\\\sqrt{-11}\cdot y_5' = (\rho^3 - \rho^8)y_1 + (\rho^5 - \rho^6)y_4 + (\rho^1 - \rho^{10})y_5 + (\rho^9 - \rho^2)y_9 + (\rho^4 - \rho^{7})y_3, \\ \sqrt{-11}\cdot y_9' = (\rho^5 - \rho^6)y_1 + (\rho^1 - \rho^{10})y_4 + (\rho^9 - \rho^2)y_5 + (\rho^4 - \rho^7)y_9 + (\rho^3 - \rho^{8})y_3,\\\sqrt{-11}\cdot y_3' = (\rho^1 - \rho^{10})y_1 + (\rho^9 - \rho^2)y_4 + (\rho^4 - \rho^7)y_5 + (\rho^3 - \rho^8)y_9 + (\rho^5 - \rho^{6})y_3, \end{cases}\end{cases}\tag{2}
\end{equation*}
the easiest combination of which includes the cyclic permutation:
\[ C) y_1' = y_4,\:y_4' = y_5, \: y_5' = y_9,\: y_9' = y_3, \: y_3' = y_1.\footnote{[Via a series of calculations one obtains the following expression for $C$ in terms of $S$ and $T$:  $C = S^6TS^2TS^6T$ (see Vol. 2 of this issue, p. 417, footnote \textsuperscript{28}).  The product of operations in the present text should always be read from left to right, that is $S^\alpha T^\beta$ will denote first performing the operation $S$ $\alpha$ times, then the operation $T$ $\beta$ times. B.-H.]} \tag{2b}\]
In fact, all 660 collineations can be given by the two symbols:
\[C^\alpha S^\beta, \quad\quad C^\alpha S^\beta T S^\gamma,\]
where $\alpha$ takes on the values $0,1,2,3,4$ and $\beta$ as well as $\gamma$ take on values $0,1,2,...,10$.\footnote{[See for instance ``\textit{Modulfunktionen}'', Vol. 2, p. 303.]}\\
This system of substitutions corresponds to a ``\textit{Problem of $y$}'' in sense of my previous essay (Matematische Annalen, Vol 15 (1879), p. 256 [=paper LVII in Vol. 2 of this issue, p. 395ff]), which goes as follows: \textit{For entire functions in $y$ which remain unchanged by the 660 collineations (2) and given by their numerical values, one seeks to compute the unknown $y$.}  In any case, this problem possesses the sought after Galois group, and I will thus occupy myself with this group for a while.
\\
\textbf{\S 4 Invariant entire functions of $y$.}
\\
It is not my intention to communicate all invariant entire functions of $y$;  in any case this would be a far reaching and possibly difficult task.  Rather I will define three of them which I will end up using later.  The first of which the function of \textit{third} degree:
\[\nabla = y_1^2y_9 + y_4^2y_3 + y_5^2y_1 + y_9^2y_4 + y_3^2y_5,\tag{3}\]
clearly the lowest invariant function.  The second is the \textit{Hessian determinant} of fifth degree:
\[ H = \left\lvert \begin{matrix} y_9 & 0 & y_5 & y_1 & 0 \\ 0 &y_3 & 0 & y_9 & y_4 \\ y_5 & 0 &y_1 & 0 & y_3 \\ y_1 & y_9 & 0 &y_4 & 0 \\ 0 & y_4 & y_3 &0&y_5\end{matrix}\right\rvert, \tag{4}\]
whose first subdeterminant will be of interest later.  The third is the \textit{function of eleventh degree}:
\[C = (y_1^{11}+y_4^{11} + y_5^{11} + y_9^{11} + y_3^{11})...,\tag{5}\]
which one can define as the numerical multiple of the sum of the eleventh powers of those values from the 660, whose $y_1$ belongs to the 660 collineations.  With respect to the first term given in $(5)$ I found that:
\[C = \frac{11}{118}(\sum y^{11}).\tag{5b}\]
From $C$ and $\nabla$ I will later construct the rational function of 33\textsuperscript{rd} degree and zeroth degree given by the following:
\[\frac{C^3}{\nabla^{11}}.\]
\\
\textbf{\S5 Eleven valued entire functions of $y$.}
\\
To find the lowest eleven-valued function of $y$ I will start by constructing a subgroup of index 11 (i.e. comprising 60 substitutions) from the collection of collineations $(2)$.  This occurs easily using the well-defined relationship between the 660 collineations and the integer coefficient linear substitutions
\[\omega' = \frac{\alpha\omega + \beta}{\gamma\omega + \delta}, \quad \quad \alpha\delta - \beta\gamma \equiv R (\text{mod } 11)\]
that differ modulo 11 and then invoking the Betti formulas that have been reproduced in paper LXXXIII, p. 79 of this issue.
\\
Clearly one can respectively assign the following two $\omega$-substitutions to the collineations $S$ and $T$:
\[\omega' = \omega + 1, \quad \quad \omega' = -\frac{1}{\omega},\]
(from which one can again construct all others through iterating and combining these two).  Then the cyclic permutations $C^\alpha$ of the five $y$:
\[(y_1,y_4,y_5,y_9,y_3)^{\alpha},\]
correspond, as one easily finds\footnote{[Using the formula for $C$ mentioned in footnote \textsuperscript{5} of this text.]}, to the iterates of
\[\omega' = 4\omega.\]
A subgroup of index 11 thus arises from the aforementioned formulas, once one connects the substitution
\[\omega' = 4\omega\]
with the following substitution of period 2:
\[\omega' = \frac{\omega - 2}{\omega - 1} = \frac{-1}{\omega - 1} + 1;\]
The subgroup as a result contains the substitution
\[\omega' = -\frac{1}{\omega}.\footnote{[For the sake of brevity, set $C^2 = U$ and $S^{-1}TS = V$.  One can thus see that $T = VC^4VCVC$; thus it becomes immediately clear that the group generated by $V$ and $C$ contains all Betti substitutions given on page 79 of the current volume under IIIa.  However, this is the extent of the substitutions it contains, since one can also take $U$ and $V$ as generator instead of $C$ and $V$, these being subject to the following relations
\[U^5 = 1, V^2 = 1, (UV)^3 = 1,\] which we know (see Dyck, Mathem, Vol 20 (1881/82), p. 35) can be used in the abstract definition of the icosahedral group. B.-H.]}\]
\textit{By returning to the $y$, we must combine the cyclic permutation of the $y$ $C$ with the collineation $S^{-1}TS$.  The resulting subgroup automatically contains the collineation $T$.}
\\
However, we already know a very simple function,that is invariant under the operations $C$ and $T$, namely the \textit{sum} of $y$:
\[p_\infty = y_1 + y_4 + y_5 + y_9 + y_3. \tag{6}\]
If one applies the collineation $S^{-1}TS$ to this, a short calculation results in
\begin{align*}
    p_0 = \frac{1}{\sqrt{-11}}\bigg\{&y_1(2(\rho^7 - \rho^1) + (\rho^9 - \rho^{10})) \\
    &+y_4(2(\rho^6 - \rho^4) + (\rho^3 - \rho^7)) \\
    &+y_5(2(\rho^2 -\rho^5) + (\rho^1 - \rho^6)) \tag{7}\\
    &+y_9(2(\rho^8-\rho^9) +(\rho^4 - \rho^2)) \\
    &+y_3(2(\rho^{10}-\rho^3) + (\rho^5 - \rho^8))\bigg\},
\end{align*}
and if one cyclically permutes the five $y$, one obtains a further four distinct expressions, which we will refer to as
\[p_1,p_2,p_3,p_4.\tag{8}\]
Since $p_\infty$ remains unchanged by 10 collineations of the subgroups, it thus is six-valued under action by collection of collineations; \textit{that is the six expressions $p$ are permuted amongst each other by the 60 Collineations of the subgroup.}  The symmetric functions of the six $p$ are invariant under all the collineations of the subgroup and are thus eleven-valued under the 660 collineations $(2)$,  they must therefore belong to the functions that are left unchanged entirely.\\
Accordingly, in order to have eleven-valued functions that are as low as possible, one calculates the smallest non-vanishing symmetric functions of $p$, i.e. the sum of squares and the sum of cubes.  One thus finds the following functions:\\
\begin{enumerate}[label = \arabic*)]
    \item \textit{The function of second degree:}
    \begin{align*}
        \varphi_0 =\quad\quad \quad \quad\quad\quad &(y_1^2 + y_4^2 + y_5^2 + y_9^2 + y_3^2) \\
        - &(y_1y_9 + y_4y_3 + y_5y_1 + y_9y_4 + y_3y_5) \tag{9}
        \\+ \frac{-1 + \sqrt{-11}}{2}&(y_1y_4 + y_4y_5 + y_5y_9 + y_9y_3 + y_3y_1),
    \end{align*}
    \item \textit{The function of third degree}
    \begin{align*}
        f_0 &= (y_1^3 + y_4^3 + y_5^3 + y_9^3 + y_3^3) \\
        &+ 3(y_1^2y_3 + y_4^2y_1 + y_5^2y_4 + y_9^2y_5 + y_3^2y_9)\\
        &-3(y_1y_4y_9 + y_4y_5y_3 + y_5y_9y_1 + y_9y_3y_4 + y_3y_1y_5) \tag{10}\\
        +\frac{1 + \sqrt{-11}}{2}&(y_1^2y_5 + y_4^2y_9+y_5^2y_3 + y_9^2y_1 + y_3^2y_4) \\
        -\frac{1+ \sqrt{-11}}{2}&(y_1y_4y_5 + y_3y_5y_9 + y_5y_9y_3 + y_9y_3y_1 + y_3y_1y_4) \\
        -(1 + \sqrt{-11})&(y_1^2y_4 + y_4^2y_5 + y_5^2y_9 + y_9^2y_3 + y_3^2y_1).
    \end{align*}
\end{enumerate}
The function $\varphi_0$ coincides with $\frac{-1 + \sqrt{-11}}{12}\sum p^2$; the function $f_0$ differs from $\frac{-\sqrt{11}\sum p^3}{6}$ on only one term, which is a numerical multiple of $\nabla$ (3).  The eleven values which $\varphi_0$ and $f_0$ take on under action of the 660 collineations, which I will call $\varphi_v$ and $f_v$ respectively, arise from $\varphi_0$ and $f_0$, when one replaces $y_{x^2}$ with $\rho^{x^2 v}\cdot y_{x^2}$ in the collineation $S^v$.  If one changes the sign of $\sqrt{-11}$ in these formulas, one obtains the expressions $\varphi_v'$ and $f_v'$ which are also eleven-valued and relate to the second collection of subgroups of index 11, which Betti also indicates.  However, since the nature of the investigations on them are structured exactly the same, as for the $\varphi_v$ and $f_v$, I will not include these in what follows.  As eleven-valued functions of zeroth degree I will later need to use both
\[\frac{f_v}{\nabla}\]
and
\[\frac{\varphi_v}{\nabla^{\frac{2}{3}}}.\]
\\
\textbf{\S6 Specialization of the $y$-problem}\\
For our special purpose we do not need to consider the \textit{generalized} problem of $y$: In the equation $J = F(z)$, which we seek out, we require only one parameter $J$ to be present.  We are thus faced with properly finding a simply extended manifold, i.e. \textit{a curve}, in the four times extended manifold of the $y_1;y_4;y_5;y_9;y_3$, which is mapped to itself by the 660 collineations, and to carry out the problem of $y$ that applies to specifically to it.
\\
Regarding this curve, we know \textit{that it must be the image of the Galois resolvent of the transformation equation.}  As I demonstrated earlier (Mathematische Annalen Vol. 14 (1878/79) [=paper LXXXII, page 55 of this volume]), we now have that the Galois resolvent is presented by a Riemann surface, which is 660 sheeted across the plane of $J$ and whose sheets are in triples at $J = 0$, in pairs at $J = 2$, and eleven-fold at $J = \infty$, and otherwise are not joined, that is their genus is $ = 26$.\footnote{[See the additional remark No. 1 on page. 166 at the end of this paper.]}  Therefore, a rational function $J$ must exist on our curve, such that it takes on every value at and only at those 660 points, which arise from the 660 collineations.  In each of these groups of 660 points each there can only be three that consist of a smaller number of multiply counted points: a group of 220 points counted three times, a group of 330 points counted two times, and a group of 60 counted eleven times.  The genus of the curve is naturally thus equal to 26.  Next, we consider the equation (1) $J = F(z)$.  Given everything else it implies that there exists a rational function of $z$ on our curve, taking on every value at 60 and only 60 points, which stem from a subgroup of collineations of index 11.  Further properties of the function $F$:  That $F(z)$ is a rational \textit{entire} function of eleventh degree, and must contain a cubic factor three times and $F(z) - 1$ must contain a biquadratic factor twice, are mere consequences of what has been established\footnote{Regarding these observations, see also the analogue observations for the transformation of seventh degree, which I developed somewhat more thoroughly in Matematische Annalen Vol. 14 (1878/79) [=paper LXXXIV, p. 113ff, 119ff of this volume].}.  That $F(z)$ is an entire function of eleventh degree follows from the fact that for the 660 points corresponding to a value of $J$ split with respect to the 60 collineations of subgroup into $11\cdot 60$, however the 60 points $J = \infty$ all arise from the collineations of the subgroup.  The other properties follow from the behavior of the 220 points of $J = 0$ and the 330 points of $J = 1$.  Among the 660 collineations, we know that there exist $2 \cdot 55$ of period 3, $55$ of period 2.\footnote{For these claims and the following see also the relevant chapter in Serret's Trait{\'e} d'alg{\`e}bre sup{\'e}rieure, vol II.[or in ``Modulfunktionen'', Vol. 1, p. 435/436.]}  Thus, each collineation of period 3 fixes $4$ points of $J = 0$, and each collineation of period 2 fixes 6 points of $J = 1$.  However, the subgroup of index 11 contains $2 \cdot 10$ collineations of period 3, 15 of period 2.  The 220 points of $J = 0$ thus separate with respect to this in $2\cdot 20 + 3 \cdot 60$ and the 330 points of $J = 1$ separate into $3 \cdot 30 + 4\cdot 60$.  And exactly this is meant by the stated properties of $F$.
\\
\textit{If one is thus successful in finding a space curve on which the function $J$ and with respect to this a function of $z$ exists in the specified way, it must directly lead to an equation $J = F(z)$ possessing all characteristic properties and thus be the equation we're looking for.}\\
Now, we have that the lowest rational function, depending only on the relations of the $y$ and unchanged by the 60 collineations of a subgroup, is given by \S5.:
\[z = \frac{f_v}{\nabla} \tag{11}\]
is a function of third degree.  As I presuppose this function should be the function taking on each value on 60 points on our curve;  I will thus introduce the hypothesis \textit{that our curve be of 20\textsuperscript{th} degree, that it lies on neither $f_v = 0$ nor on $\nabla = 0$, and does not cross the common intersection of $f_v = 0$ and $\nabla = 0$.}  I will further assume, \textit{that it does not lie on $C = 0$ and also does not cross the intersection of $C = 0, \nabla = 0$.}  Then $\frac{C^3}{\nabla^{11}}$ is a function which takes on each value at 660 associated points;  for $C = 0$ one only gets 220 separated points, and for $\nabla = 0$ only 60.  One sees: \textit{One must set $J = k\cdot \frac{C^3}{\nabla^{11}}$,} where $k$ is a numerical constant.  Now I say: \textit{If the genus of our curve is equal to 26, there directly exists only one group of multiply counted points besides the group of 220 triply counted points and the group of 60 counted eleven times, namely the 330 counted twice.}  Since if we imagine the curve as a 660-sheeted Riemann surface spread over the plane of $J$, we always have a \textit{regular} branching (Annalen XIV, pg. 458) due to the 660 transformations of the curve into itself.  And if $v$ sheets are connected anywhere, it follows that at this point all sheets are connected $v$ at a time.  One thus has
\[2p - 2 = 660\left(-2 + \sum \frac{v-1}{v}\right),\]
where the sum on the right hand side runs through the various places in the plane of $J$ where there are branchings.  Given $p = 26$, one value of $v$ equal to $3$, one equal to $11$, one can immediately show that a third $v$ shows up under the summation sign, which in any case must be equal to $2$.  Via an appropriate choice of $k$, we can also arrive at the conclusion that $J = 1$ at the relevant positions.
\\
One see how all of these hypotheses merge together.  \textit{It is the matter of finding a curve in the space of $y$ with degree 20 and genus 26, which is mapped to itself by the 660 collineations, lies on neither $f_v = 0$, $\nabla = 0$, nor finitely on $C = 0$, and does not cross the intersection of $f_v = 0$ or $\nabla = 0$, nor the intersection of $C = 0$ or $\nabla = 0$.}
\\
\textbf{\S7 The double curve of  $H = 0$}
\\
Provided our Curve of 20\textsuperscript{th} degree exists, it must lie on the surface $(4)$: $H = 0$.\footnote{I use \textit{surface} to mean every manifold, that is represented by \textit{one} equation, i.e. in the case here a three dimensional manifold, and for \textit{curve} a one dimensional manifold.}  This is true since otherwise there would be 100 points of intersection with $H = 0$, and these 100 points would have to be permuted amongst each other by the 660 collineations, which is impossible.  Now one recalls from ordinary space geometry that the Hessian of a surface of third degree has 10 nodes.  For this purpose one simultaneously sets all first subdeterminants of the Hessian determinant equal to zero.  One can proceed in the same manner in the case of 5 variables and obtains the following general theorem via use of known methods:  \textit{in the case of 5 variables the Hessian of a surface of third degree possess a double curve of 20\textsuperscript{th} degree and genus 26.} [To this end, let $H_{ik} = 0$, given in any order, be the 15 equations of fourth degree one is to set equal to zero in order to obtain the first subdeterminants of $H.$  Due to the identical relations of the $H_{ik}$ in terms of $y_1,y_4,...,y_3$ these equations are sufficiently compatible for defining a curve.  Any three of these will first determine a curve of degree $64$, from which one can remove certain components on which the remaining $H_{ik}$ fail to vanish simultaneously.  The computations proceed similarly to those carried out by Clebsch in Crelle's Journal, Vol. 59 (1874) in the computation of double points of the Hessian surface of a surface of third degree in three dimensional space.  (See also Salmon-Fiedler, \textit{Analyt. Geometrie des Raumes,} Vol. 2, 2\textsuperscript{nd} Ed. (1874), p. 331-334 and p. 535-536.)  From this one obtains that the degree equals 20. -I do not remember how I computed the genus $p = 26$ of the double curve using similarly general assumptions.  However, below on page 155 the genus of our special curve is verified as $p = 26$, and since it possesses no double points, we are presented with no reason why its genus would be smaller than that of the double curve of the Hessian surface of a general surface of third degree in four-dimensional space.  Thus one can carry over the work on the special case to that of the general case.]\footnote{[Addition upon reprinting. -see the additional remark 5 on p. 168 of this paper. K.]} \\
The surface $H = 0$ will also possess such a double curve, provided special relations don't affect this.  This double curve will be mapped to itself by the 660 collineations, since the same holds true of the surface $H = 0$.  Can one doubt that the double curve is the curve we seek?  To this end there only remain two more things to show:  first of all that the numbers 20 and 26 for degree and genus in the general case are not subject to any modifications in this special case, and secondly, that our curve also possess the other, negative properties that we have stated.
\\
To anticipate the proof that I will get to immediately, I will already state the theorem here:
\\
\textit{The curve of 20\textsuperscript{th} degree we seek is the double curve of the Hessian surface $H = 0.$}
\\
For the proof we will first construct all subdeterminants of $H$ and set them equal to zero.  One thus gets a system of equations, which I will refer to as 
\[H_{ik} = 0 \tag{12}\]
and whose equations arise from the following three
\[\tag{13}\begin{cases}0 = y_4y_5y_9y_3 - y_1^2y_5y_3 + y_1^2y_4^2 + y_3^3y_1, \\ 0 = y_1^2y_5y_9 - y_4^2y_5y_3 - y_3^2y_1y_9,\\
0 = y_4^3y_9 + y_9^3y_5 + y_3^3y_1\end{cases}\]
up to cyclic permutation of the $y$.
\\
I will now designate the five points, where four of the five $y$ vanish, as the points I,IV,V,IX,III.  One then sees immediately that:
\\
\textit{The five points } I, IV, V, IX, III \textit{belong to our curve.} And since they clearly do not belong to the surfaces $(5)$: $C = 0$ and $f_v = 0$ (\S5), it follows:
\\
\textit{Our curve lies on neither $C = 0$ nor $f_v = 0$.}
\\
If one sets one of the $y$ equal to zero in $(13)$ and the remaining equations $H_{ik} = 0$, it follows that any three more of the remaining $y$ must also vanish.  Therefore:
\\
\textit{Every plane $y_{k^2} = 0$ intersects our curve at only four points, namely at those four of the five points} I, IV, V, IX, III, \textit{which are not named after the index of $y_{k^2}$.}
\\
Next, one takes one of the five points, say the point III, and one introduces a series expansion there, by setting $y_3 = 1$, $y_5 = dt$.  The equations\footnote{[For this calculation it even suffices to only use those equations that arise from cyclic permutations of the indices of from the last of the three equations in (13). B-H.]} $H_{ik} = 0$ then give us the expansion up to terms of tenth degree:
\[y_1 = dt^{10}, \quad y_4 = dt^6, \quad y_5 = dt, \quad y_9 = -dt^3, \quad y_3 = 1\]
which leads to the following table for the relations between the $y$ at all five points:
\[\begin{tabularx}{\textwidth}{ X|X|X|X|X|X }
     & $y_1$ & $y_4$ & $y_5$ & $y_9$ & $y_3$ \\
    \hline 
    I & 1 & $dt^{10}$ & $dt^6$ & $dt$ & $-dt^3$  \\
    \hline
    IV & $-dt^3$ & 1 & $dt^{10}$ & $dt^6$ & $dt$
    \\
    \hline
    V & $dt$ & $-dt^3$ & 1 & $dt^{10}$ & $dt^6$
    \\
    \hline
    IX & $dt^6$ & $dt$ &$-dt^3$ & 1 & $dt^{10}$
    \\
    \hline
    III & $dt^{10}$ & $dt^6$ & $dt$ & $-dt^3$ & 1
    \\
    \end{tabularx} \tag{14}\]
One sees:
\\
\textit{The five points are simple points of our curve.}\\
However especially:
\\
\textit{The curve is of degree 20.}
\\
Since the sum of the exponents of $dt$ for a single $y$ appearing in the table is $3 + 1 + 6 + 10 = 20.$\\
Each of our five points remains invariant under the collineations $S: y'_{k^2} = \rho^{k^2}y_{k^2}$ and are permuted amongst each other by the cyclic permutation $C$.  Thus our 5 points only give rise to 60 under the 660 collineations.  $\nabla$ vanishes at all of these points because it vanishes at I for instance.  However, $\nabla$ does not vanish identically.  Since if we plug the suitable series expansion (say at I) in (14) to $\nabla$, we get (given our chosen measure of accuracy) $\nabla = dt$.  Thus:
\\
\textit{$\nabla = 0$ has precisely 60 points of intersection with our curve, and neither $C$ nor $f_v$ vanish at these points.}
\\
To have all of the characteristic properties of the curve we're seeking it only remains to show, \textit{that the genus is equal to 26.}  This holds relatively easily by examining certain inequalities.  [By viewing our double curve of degree 20 as being contained in the intersection of three surfaces of fourth degree, which we may call $H_1 =0$, $H_2 = 0, H_3 = 0$, we know we can find:]\footnote{[In the reprint this replaced an inaccurate formulation found in the original. -One can use the following approach: Let $u_y = 0, v_y = 0$ be any two planes, $(u,v,H_1,H_2,H_3)$ the functional determinant of the given forms, one can obtain the following everywhere finite integrals belonging to our curve of 20\textsuperscript{th} degree:
\[\int \Phi_7\frac{u_ydv_y - v_ydu_y}{(u,v,H_1,H_2,H_3)}\]
where we interpret $\Phi_7$ as a form of seventh degree in the variables $y$.  The surface $\Phi_7 = 0$ is thus subject to the condition that it intersects the curve of 20\textsuperscript{th} degree in those $k$ rigid points, in which the curve encounters the remaining components of the intersection curve of $H_1 = 0, H_2 = 0, H_3 = 0$.  In addition $\Phi_7 = 0$ intersects the curve of degree 20 in $2p-2$ movable points.  From this we obtain the relation
\[k+2p - 2 = 7\cdot 20 = 140,\]
thus in particular $p \leq 71.$  Furthermore, fixing $p = 26$, it follows that our $\Phi_7 = 0$ must intersect the curve of degree 20 in $140 - (2\cdot 26 - 2) = 90$ rigid points, if we would like to have the integral of the first kind as above exist. K.]}
\[p < 71\]
On the other, we know that $p$ is no smaller than the genus of a curve that decomposes into straight lines, that is:
\[p \geq -19.\]
Finally it is possible to set up the following equation, as showed in the previous section:
\[2p - 2 = 660 \left(-2 + \sum \frac{v - 1}{v}\right).\]
In this one equation one $v$ must be taken as 3 due to the intersection with $C = 0$, and another $v$ equal to $11$ due to the intersection with $\nabla = 0$.  If there are no other $v$, then it would follow that
\[2p - 2 = - 280, \quad \quad p = -139.\]
However, if one were to take the third $v$ as greater than or equal to $3$, we'd get that
\[2p - 2 \geq 160, \quad \quad p \geq 81.\]
\textit{That means there only exists one more $v$ which is equal to 2; meaning $p = 26$}, Q.E.D.
\\
\textbf{\S8 The equation $J = F(z).$}
\\
What remains is a simple calculation in degree to set up the equation $J = F(z)$.  We can set
\[J = k\cdot \frac{C^3}{\nabla^{11}}, \quad \quad z_v = \frac{f_v}{\nabla}\tag{15}\]
and first consider as having undetermined coefficients\footnote{The $k$ on the right-hand side of the equation is in fact the same as in the $k$ in (15), since in $\frac{C^3}{\nabla^{11}}$ the numerator has the term $y_1^{33}$ with 1 which appears multiplied, just like in $z^{11}$.}:
\[J = k(z^2 + Az + B)(z^3 + az^2 + bz + c)^3 \tag{16}\]
or:
\[J - 1 = k(z^3 + Az^2 + Bz + \Gamma)(z^4 + \alpha z^3 + \beta z^2 + \gamma z + \delta)^2. \tag{16b}\]
Next we will plug in one of the series expansions in (14) into $C$, $\nabla$ and $f_v$.  Thus, except for terms of higher than tenth degree, we have that:
\[\begin{cases} C = 1, \quad \quad \nabla = dt \\ f_v = \rho^{9v} + \frac{1 + \sqrt{-11}}{2}\rho^{2v}\cdot dt^2 -2\rho^{4v}\cdot dt^3 + \frac{1 + \sqrt{-11}}{2}\rho^{6v}\cdot dt^4 \\ \quad\quad +(1 + \sqrt{-11})\cdot dt^5 - \frac{1 + \sqrt{-11}}{2}\cdot \rho^{10v}\cdot dt^6 + 3\rho^{3v}\cdot dt^8 \\ \quad \quad +2\rho^{5v}\cdot dt^9 - (1 + \sqrt{-11})\rho^{7v}\cdot dt^{10}.
\end{cases}\tag{17}\]
Up to higher order terms, plugging this into (16) and (16b) gives us the following:
\[A = - 3,\quad B = 5 - \sqrt{-11}, \quad a = 1, \quad b = -3\frac{1 + \sqrt{-11}}{2}, \quad c = \frac{7 - \sqrt{-11}}{2};\tag{18}\]
\[A = 4, \quad  B = \frac{7 - 5\sqrt{-11}}{2}, \quad \Gamma = 4 - 6\sqrt{-11},\quad \alpha = -2, \quad \beta = 3\cdot \frac{1 - \sqrt{-11}}{2}, \quad \gamma = 5 + \sqrt{-11}, \quad \delta = -3 \cdot \frac{5 + \sqrt{-11}}{2}, \tag{18b}\]
and the values we get for $J$ and $J-1$ from this do indeed coincide (which in turn is the source of plenty of verifications) if we set
\[k = -\frac{1}{1728}. \tag{19}\]
\textit{Therefore the finished equation $J = F(z)$ is given by the following:}
\begin{align*}
    J: J - 1: 1 &= (z^2 - 3z + (5 - \sqrt{-11}))\cdot\\
    &\cdot \left(z^3 + z^2 - 3\cdot\frac{1 + \sqrt{-11}}{2}\cdot z + \frac{7 - \sqrt{-11}}{2}\right)^3 \\
    &:\left(z^3 + 4z^2 + \frac{7-5\sqrt{-11}}{2}\cdot z + (4 - 6\sqrt{-11})\right)\cdot \tag{20} \\
    &\cdot \left(z^4 -2z^3 + 3\cdot\frac{1- \sqrt{-11}}{2}\cdot z^2 + (5 + \sqrt{-11})z - 3\cdot\frac{5 + \sqrt{-11}}{2}\right)^2 \\
    &: -1728.
\end{align*}
\\
\textbf{\S9 The second form of the equation of eleventh degree}\\
Besides the equation we just obtained we can obtain a second with a similar amount of ease, given one starts with the function of \textit{second} degree $\varphi_v$ (9) instead of the eleven-valued function of \textit{third} degree $\varphi_v$ (9).  I will first prove the following theorem: \\
\textit{Subjected to the relations $H_{ik} = 0$ (12), we can reduce the collection of entire functions of $y$ invariant under the 660 collineations to entire functions of $\nabla$ and $C$.}
\\
Setting an invariant function of $y$ equal to zero gives us a surface which either contains our curve - and thus would be identically zero modulo the $H_{ik} = 0$, or intersects the curve at points that are permuted among one another under the 660 collineations.  Among these we can find the 60 points $\nabla = 0$ a certain number of times, as well as the 220 points $C = 0$ arbitrarily often, and any groups of 660 separated points given by $C^3 - \lambda\nabla^{11} = 0$, where $\lambda$ is a suitable constant.  On the other hand the group of 330 points $J = 1$ can only show up an \textit{even} number of times, because the total number of points must be divisible by 20, but 330 is only divisible by 10.  However, this group is represented by a combination of a $C$ and $\nabla$ doubly counted, since $J = 1$ is through $C^3 + 1728\nabla^{11} = 0;$  all points of intersection can thus be cut out with the right multiplicity but setting a suitable entire function of $\nabla$ and $C$ equal to zero, Q.E.D.
\\
\textit{Therefore thanks to the relations $H_{ik} = 0$ the eleven $\varphi_v$ satisfy an equation of eleventh degree, whose coefficients are given by entire functions of $\nabla$ and $C$.}
\\
With respect to the degree of the functions to be considered we can ascribe undetermined numerical factors:
\[\varphi^{11} \alpha\nabla^2\cdot \varphi^8 + \beta\nabla^4\cdot \varphi^5 + \gamma\nabla C\cdot \varphi^4 + \delta \nabla^6\cdot \varphi^2 + \epsilon\nabla^3 C\cdot \varphi + \zeta\cdot C^2 = 0.\tag{21}\]
The numerical factors are once again determined with the help of the series expansions (14).  Using them one has:
\[\varphi_v = \rho^{6v} - \rho^{8v}\cdot dt + \rho^{10v}\cdot dt^2 + \frac{1 - \sqrt{-11}}{2}\rho^v\cdot dt^3 + \frac{1-\sqrt{-11}}{2}\rho^{3v}\cdot dt^4 + ...\tag{22}\]
and thus it follows that:
\[\alpha = -22, \beta = 11(9 - 2\sqrt{-11}), \gamma = 11, \delta = 88\sqrt{-11}, \epsilon = \frac{11(-3 + \sqrt{-11})}{2}, \zeta = -1.\tag{23}\]
I also to set 
\[\frac{\varphi_v}{\nabla^{\frac{2}{3}}} = \xi_v \tag{24}\]
and to introduce for $\frac{C^3}{\nabla^{11}}$: $-1728 J = -1728\frac{g_2^3}{\Delta}$. \textit{Thus the new equation of eleventh degree takes on the following form:}
\[\xi^{11} -22\cdot \xi^8 + 11(9 - 2\sqrt{-11})\xi^5 - 11 \cdot\frac{11g_2}{\sqrt[3]{\Delta}}\cdot\xi^4 + 88\sqrt{-11}\cdot\xi^2 - 11\cdot \frac{-3 + \sqrt{-11}}{2}\cdot \frac{12g_2}{\sqrt[3]{\Delta}}\cdot \xi - \frac{144g_2^2}{\sqrt[3]{\Delta^2}} = 0. \tag{25}\]
I remains only to show how this equation is connected with equation (20).  The surface $\varphi_v = 0$ intersects with 40 points of our curve which are permuted under one another by the 60 collineations of the subgroup.  Due to the above this can only be those $2\cdot 20$ points, which are each fixed by 2 collineations of period 3, that is the same points, for which equation (20) gives:
\[z^2 - 3z + (5 - \sqrt{-11}) = 0.\]
In fact, the series expansions (14) show that the following relation holds (naturally still given the $H_{ik} = 0$:
\[\varphi_v^3 = f_v^2 - 3f_v\nabla + (5 - \sqrt{-11}).\tag{26}\]
and \textit{that one can get equation (20) from equation equation (25) by setting:}
\[\xi^3 = z^2 - 3z + (5 - \sqrt{-11}).\tag{27}\]
The direct verification of this claim, which once again comprises a series of estimates of the numerical coefficients, presents no difficulties.
\\
\textbf{\S10. Summary of the previous results.}
\\
Summarizing, we have obtained the following results for \textit{the transformation of eleventh degree of elliptic functions:}
\begin{enumerate}[label = \arabic*)]
    \item \textit{The Galois resolvent of 660\textsuperscript{th} degree can be written in the following way:}  One takes the five related quantities
    \[y_1:y_4:y_5:y_9:y_3\]
    subject to the 15 relations $H_{ik} = 0$ (see (12) and (13)) and set:\footnote{If one wishes to examine the $y$ themselves rather than relations of them, one could write $C = 12g_2,$ $\nabla = - \sqrt[11]{\Delta}$, that is one would have to adjoin $g_2$ and $\sqrt[11]{\Delta}$. }
    \[\frac{-C^2}{1728\nabla^{11}} = J,\]
    where $\nabla$ denotes the function of third degree (2), and $C$ denotes the function of eleventh degree (4).  If one has found \textit{a} solution system of these equations, the remaining ones can be found using the collineations in \S3.
    \item \textit{There exist two simplest forms of the resolvent of eleventh degree.}  The one which we first considered alone, is given by (20):
    \begin{align*}
    J: J - 1: 1 &= (z^2 - 3z + (5 - \sqrt{-11}))\cdot\\
    &\cdot \left(z^3 + z^2 - 3\cdot\frac{1 + \sqrt{-11}}{2}\cdot z + \frac{7 - \sqrt{-11}}{2}\right)^3 \\
    &:\left(z^3 + 4z^2 + \frac{7-5\sqrt{-11}}{2}\cdot z + (4 - 6\sqrt{-11})\right)\cdot \\
    &\cdot \left(z^4 -2z^3 + 3\cdot\frac{1- \sqrt{-11}}{2}\cdot z^2 + (5 + \sqrt{-11})z - 3\cdot\frac{5 + \sqrt{-11}}{2}\right)^2 \\
    &: -1728;
\end{align*}
its 11 roots are given by the formula:
\[z_v = \frac{f_v}{\nabla},\]
where $f_v$ is defined in equation (10).
\\
The second form is introduced by (25):
\[0 = \xi^{11} -22\cdot \xi^8 + 11(9 - 2\sqrt{-11})\xi^5 - 11 \cdot\frac{11g_2}{\sqrt[3]{\Delta}}\cdot\xi^4 + 88\sqrt{-11}\cdot\xi^2 - 11\cdot \frac{-3 + \sqrt{-11}}{2}\cdot \frac{12g_2}{\sqrt[3]{\Delta}}\cdot \xi - \frac{144g_2^2}{\sqrt[3]{\Delta^2}};\]
and its roots are:
\[\xi_v = \frac{\varphi_v}{\nabla^{\frac{2}{3}}},\]
where $\varphi_v$ are understood as the functions (9).
\end{enumerate}
\textbf{\S 11 Connection with equation of twelfth degree.}
\\
I now wish to still show, how the quantities $y$ relate to the multiplier equation of twelfth degree, which I recently wrote\footnote{See Mathematische Annalen, Vol. 15 (1878/79) [=preceding printed Note No. LXXXV, p. 139].} about, and additionally how one can accordingly solve the aforementioned problem of 660\textsuperscript{th} degree, i.e. the equations of eleventh degree, in a transcendental way.  First of all I will include the calculated equation of twelfth degree here:
\[z^{12} - 90\cdot11\cdot\sqrt[2]{\Delta}\cdot z^6 + 40\cdot11\cdot12g_2\cdot\sqrt[3]{\Delta}\cdot z^4 - 15\cdot11\cdot216g_3\cdot\sqrt[4]{\Delta}\cdot z^3 + 2\cdot11\cdot(12g_2)^2\cdot\sqrt[6]{\Delta}\cdot z^2 - 12g_2\cdot216g_3\cdot\sqrt[12]{\Delta}\cdot z - 11\cdot \Delta = 0,\]
and will first and foremost explain how its roots can be represented as functions of the periodic relation $\frac{\omega_1}{\omega_2} = \omega$ of the elliptic integral, i.e. as functions of $q = e^{i\pi\omega}$.  For (28) a Jacobian equation is known.  If one accordingly sets the following:
\[\begin{cases}\sqrt{z_\infty} = \sqrt{-11}\cdot A_0, \\ \sqrt{z_v} = A_0 + \rho^v A_1 + \rho^{4v}A_4 + \rho^{5v}A_5 + \rho^{9v}A_9 + \rho^{3v}A_3, \quad \quad (v = 0,1,...,10)\end{cases}\tag{29}\]
(I only deviate from the Jacobian description by choosing the indices for $A$ to be quadratic residues modulo 11), one obtains the following formulas for \textit{a} value system of the $A$ using familiar methods.
\[\begin{cases}\mu A_0 = q^{\frac{121}{132}}\cdot \sum\limits_{-\infty}^{+\infty}(-1)^{h+1}\cdot q^{33h^2 + 55h + 22}, \\ \mu A_1 = q^{\frac{1}{132}}\cdot\left[\sum\limits_{-\infty}^{+\infty}(-1)^h \cdot q^{33h^2 +h} + \sum\limits_{-\infty}^{+\infty}(-1)^{h+1}\cdot q^{33h^2 + 13h + 14}\right] \\ \mu A_4 = q^{\frac{37}{132}}\cdot\left[\sum\limits_{-\infty}^{+\infty}(-1)^h \cdot q^{33h^2 +13h + 1} + \sum\limits_{-\infty}^{+\infty}(-1)^{h+1}\cdot q^{33h^2 + 31h + 7}\right] \\ \mu A_5 = q^{\frac{49}{132}}\cdot\left[\sum\limits_{-\infty}^{+\infty}(-1)^h \cdot q^{33h^2 + 37h + 10} + \sum\limits_{-\infty}^{+\infty}(-1)^{h+1}\cdot q^{33h^2 + 7h}\right] \\ \mu A_9 = q^{\frac{97}{132}}\cdot\left[\sum\limits_{-\infty}^{+\infty}(-1)^{h+1} \cdot q^{33h^2+19h + 2} + \sum\limits_{-\infty}^{+\infty}(-1)^{h}\cdot q^{33h^2 + 25h + 4}\right] \\ \mu A_3 = q^{\frac{25}{132}}\cdot\left[\sum\limits_{-\infty}^{+\infty}(-1)^h \cdot q^{33h^2+49h + 18} + \sum\limits_{-\infty}^{+\infty}(-1)^{h}\cdot q^{33h^2 + 61h + 28}\right]\end{cases}\tag{30}\]
where $\mu$ denotes the factor of proportionality $\sqrt{\frac{\omega_2}{\pi}}$.\footnote{[So that the formulae (29), (30) can exist together, one must put the roots $z_r$ of the multiplier equation in the correct order, which deviates from the ordinary one in following sense:
\[z_\infty = \left(\frac{2\pi}{\omega_2}\right)\cdot(-11)\cdot q^{\frac{11}{6}}\cdot\Pi(1-q^{22\lambda})^2\]
\[z_r = \left(\frac{2\pi}{\omega_2}\right)\cdot\rho^{2\nu}\cdot q^{\frac{1}{66}}\cdot\Pi(1-\rho^{2\nu\lambda}q^{\frac{2}{11}\lambda})^2\]
\[(\nu = 0,1,...,10).\]
The general expression for the series expansion (30) with non-zero index $k^2$ (mod $11$) is as follows:
\[\mu\cdot A_{k^2} = q^{\frac{1}{132}}\cdot\left[\sum\limits_{-\infty}^{+\infty}(-1)^h\cdot q^{33h^2 - (11+12k)h + \frac{2(k+1)(6k+5)}{11}} + \sum\limits_{-\infty}^{+\infty}(-1)^h\cdot q^{33h^2 - (11-12k)h + \frac{2(k-1)(6k-5)}{11}} \right]\]
In the text the summation letter $h$ is replaced by $h+1$ or $h+2$ only a few times. B.-H.]}
\\
There are $660\cdot 24$ such value systems, 660 due to the Galois group of the modular equation, and 24 due to twelfth root appearing in (28) and the square root appearing in (29).  But one also easily sees that 24 value systems at a time differ solely by a 24\textsuperscript{th} root of unity, \textit{meaning that the relations on the $A$ are simply 660-valued}.  In fact, if one lets $\omega$ in the formulas (30) grow by 11 units, all $A$ get a 24\textsuperscript{th} root of unity, although this is the same one for all $A$, and only gives rise to the common factor $q^{\frac{1}{132}}$ on the right hand side. Similarly, from the discussion above, we had that the relations on the $y$ were 660-valued.  This gives rise to the following possibility: \textit{rationally expressing the relations on the $y$ through the relations on the $A$ and vice-versa,} and this is the precise formulation of the problem with which we will occupy ourselves with now.
\\
\textbf{\S 12 Related Formulas.}
\\
I will summarize the observations required for deriving this exceedingly easy result by a series of individual remarks.  In doing so I want to explicitly emphasize that I would barely have been led along this train of thought if there wasn't a very similar result for the transformation degrees of 5 and 7.\footnote{In the future, I intend to show how the analogous theorem comes about for arbitrary degrees. [See also Mathematische Annalen., Vol. 17 (1881) [=paper LXXXIX in the present volume, p. 190, footnote \textsuperscript{11}.]]}
\begin{enumerate}[label = \arabic*)]
    \item The 660 value systems of the $A_0:A_1:A_4:A_5:A_9:A_3$ arise according to my representation in Vol. 15 of Mathematische Annalen (1879), p. 276 [=paper LVII, Vol. 2 of this issue, p. 417] from the 660 already known collineations, from which I will take the following two:
    \[\begin{cases}S' \; A_0' = A_0, \; A_1 = \rho A_1, \; A_4' = \rho^4A_4, \; A_5' = \rho^5A_5, \; A_9' = \rho^9A_9, \; A_3' = \rho^3A_3, \\ C' \; A_0' = A_0, \; A_1 = \rho A_4, \; A_4' = A_5, \; A_5' = A_9, \; A_9' = A_3, \; A_3' = A_1. \end{cases}\tag{31}\]
    \item A well-defined correspondance between the relations on the $A$ and the relations on the $y$ can be achieved in 660 different ways.  Since if one establishes one given correspondence, one can naturally act on either the $y$ or the $A$ with any one of the 660 collineations.  Therefore one can conclude that the iterations of $S:$
    \[y_1' = \rho y_1, y_4' = \rho^4y_4, y_5' = \rho^5y_5, y_9' = \rho^9y_9, y_3' = \rho^3y_3\]
    correspond to the iterations of $S'$ and the cyclic permutations
    \[(y_1y_4y_5y_9y_3)\]
    correspond to the iterations of $C'$.
    \item I now claim, under this condition, that the cyclic substitution C:
    \[y_1' = y_4, y_4' = y_5, y_5' = y_9, y_9' = y_3, y_3' = y_1\]
    corresponds necessarily to the cyclic permutation $C'$ itself, and not one of its iterates.  Since if:
    \[\frac{A_1}{A_0} = R(y_1,y_4,y_5,y_9,y_3),\]
    where $R$ represents a rational function of zeroth degree.  If we multiply each $y_{k^2}$ by $\rho^{k^2}$, according to the collineation $S$, then $R$ picks up a eleventh root of unity $\rho^v$ as a factor.  By using the collineation $C$, we get:
    \[R(y_4,y_5,y_9,y_3,y_1).\]
    If we write $\rho^{k^2}y_{k^2}$ instead of $y_{k^2}$, $\rho^{4v}$ must show up as a factor.  Therefore $R(y_4,y_5,y_9,y_3,y_1)$ must be equal to $\frac{A_4}{A_0}$, QED.
    \item The five points (I,IV,V,IV,III) (see \S7.) on the $y$-curve were characterized by the fact that they were simultaneously invariant under the collineation $S$ (and its iterates).  On the curve of the $A$ (if this geometric terminology is permitted), they correspond to five points which permit the substitution $S'$.  \textit{Clearly these are exactly those five points where $A_0$ and four of the remaining $A$s vanish.}  First of all, this is due to the obvious fact fact that these points remain unchanged under the collineation $S'$, and secondly they belong to the curve of the $A$s.  If one first takes the common factor $q^{\frac{1}{132}}$ on the right hand side of the expressions in (30) and sets $q = 0$, one obtains
    \[A_0 = 0, A_1 \neq 0, A_4 = A_5 = A_9 = A_3 = 0, \]
    so that our claim holds for one of the five points; however, we obtain the remaining four points using the cyclic permutation $C'$.  I will call these points I',IV',V',IX',III'.
    \item Due to remarks $2)$ and $3)$ we can assign any point I',...,III' to point I; if for instance one associates the point IV' to I, we necessarily have that IV,V,IX,III respectively correspond to V',IX',III',I'.
    \item  $A_0$ can only vanish at the points I',...,III'.  Since if $A_0$ is equal to zero, it follows from equation (29) that one of the roots $z$ is equal to zero, i.e., by (28), $\Delta = 0$ or $J = \infty$.  There exist 60 points $J = \infty$, but for only five of these is it possible for the single root $z$ to be zero.  Since at the 5 points I'...III', thanks to (29) only $z_{\infty}$ vanishes, and none of the other roots $z_v$.
    \item After I removed the common factor $q^{\frac{1}{132}}$ on the right hand side of (30), I wish to set $q^{\frac{2}{11}} = -ds$. We first have that
    \[A_0, A_1,A_4,A_5,A_9,A_3\]
    are approximately proportional to:
    \[-ds^5, 1, -ds^7, -ds^2,ds^{15},ds\]
    respectively.
    This gives the following table for the behavior of the $A$ at the points I',...,III':
    \[\begin{tabularx}{\textwidth}{ X|X|X|X|X|X|X }
     & $A_0$ & $A_1$ & $A_4$ & $A_5$ & $A_9$ & $A_3$ \\
    \hline 
    I' & $-ds^5$ & $1$ & $-ds^7$ & $-ds^2$ & $+ds^{15}$ & $+ds$  \\
    \hline
    IV' & $-ds^5$ & $+ds$ & $1$ & $-ds^7$ & $-ds^2$ & $+ds^{15}$
    \\
    \hline
    V' & $-ds^5$ & $+ds^{15}$ & $+ds$ & 1 & $-ds^7$ & $-ds^2$
    \\
    \hline
    IX' & $-ds^5$ & $-ds^2$ &$+ds^{15}$ & 1 & $-ds^7$
    \\
    \hline
    III' & $-ds^5$ & $-ds^7$ & $-ds^2$ & $+ds^{15}$ & 1
    \\
    \end{tabularx} \tag{32}\]
    \item \textit{The curve of the $A$ has 25\textsuperscript{th} degree.}  This is because the sum of the exponents of the $ds$ in the column corresponding to $A_0$ is 25 in the table above.  However, the sum of exponents is also 25 for any other column belonging to a different $A$.  \textit{Thus none of the other $A$ vanish away from the points I'...III'}.\\
    In particular, it follows that
    \[A_0^5 + A_1A_4A_5A_9A_3 = 0, \tag{33}\]
    which is a relation which Brioschi makes occasional use of.\footnote{Sopra una classe di equazioni modulari.  Annali di Matematic, ser. 2, t. IX, pg. 167 ff. [=Opere matematiche, No. LXXV, part II. p. 193] [Equation (33) arises from the following observation: If the curve of $A$ did not lie on the surface represented by the $A$, it would have 125 points of intersection in common with it.  The total order of vanishing of the form $A_0^5 + A_1A_4A_5A_9A_3$ on the curve of $A$ is greater than 125 due to our formulas (25) and this curve is irreducible. K.]}
    \item Next one considers the following relations on the $y$:
    \[\frac{y_4}{y_5},\frac{y_5}{y_9},\frac{y_9}{y_3},\frac{y_3}{y_1},\frac{y_1}{y_4}.\]
    Since the $y$ do not vanish away from the points I...III, one obtains the following table for when these functions are zero and infinity (see (14)):
    \[\begin{tabularx}{\textwidth}{ X|X|X|X|X|X }
    I & $+dt^{4}$ & $+dt^{5}$ & $-dt^{-2}$ & $-dt^{3}$ & $+dt^{-10}$  \\
    \hline
    IV & $dt^{-10}$ & $+dt^{4}$ & $+dt^{5}$ & $-dt^{-2}$ & $-dt^{3}$
    \\
    \hline
    V & $-dt^{3}$ & $+dt^{-10}$ & $+dt^{4}$ & $+dt^{5}$ & $-dt^{-2}$
    \\
    \hline
    IX & $-dt^{-2}$ & $-dt^{3}$ &$+dt^{-10}$ & $+dt^{4}$ & $+dt^{5}$
    \\
    \hline
    III & $+dt^{5}$ & $-dt^{-2}$ & $-dt^{3}$ & $+dt^{-10}$ & $+dt^{4}$
    \\
    \end{tabularx} \tag{34}\]
    \item On the other hand one can consider the following relations of the $A$:
    \[-\frac{A_0}{A_1}, -\frac{A_0}{A_4}, -\frac{A_0}{A_5},-\frac{A_0}{A_9}, -\frac{A_0}{A_3}.\]
    The zeroes and infinities are similarly only at the points I'...III', and in particular one obtains the following table from (32):
    \[\begin{tabularx}{\textwidth}{ X|X|X|X|X|X }
    I' & $+ds^{5}$ & $-ds^{-2}$ & $-ds^{3}$ & $+ds^{-10}$ & $+ds^{4}$  \\
    \hline
    IV' & $ds^{4}$ & $+ds^{5}$ & $-ds^{-2}$ & $-ds^{3}$ & $+ds^{-10}$
    \\
    \hline
    V' & $+ds^{-10}$ & $+ds^{4}$ & $+ds^{5}$ & $-ds^{-2}$ & $-ds^{3}$
    \\
    \hline
    IX' & $-ds^{3}$ & $+ds^{-10}$ &$+ds^{4}$ & $+dt^{5}$ & $-ds^{-2}$
    \\
    \hline
    III' & $-ds^{-2}$ & $-ds^{3}$ & $+ds^{-10}$ & $+dt^{4}$ & $+dt^{5}$
    \\
    \end{tabularx} \tag{35}\]
    \item Now one assigns IV' to I, and thus, by 5), V',IX',III',I' to IV,V,IX,III respectively.  Then, by a glance at the table, we have that
    \[\frac{y_4}{y_5},\frac{y_5}{y_9},\frac{y_9}{y_3},\frac{y_3}{y_1},\frac{y_1}{y_4}.\]
    and
    \[-\frac{A_0}{A_1}, -\frac{A_0}{A_4}, -\frac{A_0}{A_5},-\frac{A_0}{A_9}, -\frac{A_0}{A_3}.\]
    have the exact same number and locations of zeroes and poles, and can therefore be set equal to each other.  \textit{One thus has the following formulas, which complete the posited problem and supplement the results summarized in \S10 in the sense that was indicated:}
    \[\frac{y_4}{y_5} = -\frac{A_0}{A_1}, \frac{y_5}{y_9} = -\frac{A_0}{A_4}, \frac{y_9}{y_3} = -\frac{A_0}{A_5}, \frac{y_3}{y_1} = -\frac{A_0}{A_9}, \frac{y_1}{y_4} = -\frac{A_0}{A_3}.\tag{36}\]
\end{enumerate}
Ebenhausen, the 15. August 1879.
 \end{document}